\documentclass[12pt,a4paper]{amsart}

\DeclareMathAlphabet{\mathpzc}{OT1}{pzc}{m}{it}
\usepackage{amsfonts,amsmath,amssymb,indentfirst,mathrsfs,amscd}
\usepackage[mathscr]{eucal}

\date{}
\author{Marina Logares}
\address{Departamento de Matem\'aticas \\
  CSIC \\ Serrano 121
  \\ 28006 Madrid \\ Spain}

\email{marina.logares@mat.csic.es}
\title{Betti numbers of parabolic $U(2,1)$-Higgs bundles moduli spaces.}

\thanks{Partially supported by Ministerio de Educaci\'{o}n y Tecnolog\'{\i}a
through Acci\'{o}n Integrada Hispano-Lusa HP-2000-0015 and by The
European Contract Human Potential Programme, Research Training
Network HPRN-CT-2000-00101.}



\DeclareMathOperator{\rk}{rk \,}          
\DeclareMathOperator{\Jac}{Jac\,}        
\DeclareMathOperator{\di}{div\,}           
\DeclareMathOperator{\pdeg}{pardeg \,}      
\DeclareMathOperator{\pmu}{par\mu \,}       
\DeclareMathOperator{\PH}{ParHom \,}        
\DeclareMathOperator{\SPH}{SParHom\,}       
\DeclareMathOperator{\SPE}{SParEnd\,}       
\DeclareMathOperator{\Hom}{Hom\,}           
\DeclareMathOperator{\id}{Id\,}             

\DeclareMathOperator{\Coeff}{Coeff\,}

\begin{document}

\newtheorem{thm}{Theorem}
\newtheorem{prop}[thm]{Proposition}
\newtheorem{lem}[thm]{Lemma}
\newtheorem{cor}[thm]{Corollary}

\theoremstyle{definition}
\newtheorem{defn}[thm]{Definition}
\newtheorem{ex}[thm]{Example}
\newtheorem{as}[thm]{Assumption}

\theoremstyle{remark}
\newtheorem{rmk}[thm]{Remark}

\theoremstyle{remark}
\newtheorem*{prf}{Proof}

\newcommand{\iacute}{\'{\i}} 
\newcommand{\norm}[1]{\lVert#1\rVert} 

\newcommand{\cC}{\mathcal{C}}
\newcommand{\cD}{\mathcal{D}}
\newcommand{\cG}{\mathcal{G}} 
\newcommand{\cO}{\mathcal{O}} 
\newcommand{\cM}{\mathcal{M}} 
\newcommand{\cN}{\mathcal{N}} 
\newcommand{\cP}{\mathcal{P}} 
\newcommand{\cS}{\mathcal{S}} 
\newcommand{\cU}{\mathcal{U}} 
\newcommand{\cX}{\mathcal{X}}

\newcommand{\mM}{\mathscr{M}} 

\newcommand{\CC}{\mathbb{C}} 
\newcommand{\HH}{\mathbb{H}} 
\newcommand{\RR}{\mathbb{R}} 
\newcommand{\ZZ}{\mathbb{Z}} 

\renewcommand{\lg}{\mathfrak{g}} 
\newcommand{\lh}{\mathfrak{h}} 
\newcommand{\lu}{\mathfrak{u}} 
\newcommand{\la}{\mathfrak{a}} 
\newcommand{\lb}{\mathfrak{b}} 
\newcommand{\lm}{\mathfrak{m}} 
\newcommand{\lgl}{\mathfrak{gl}} 
\newcommand{\ext}{\mathrm{ext}\,} 

\newcommand{\imat}{\sqrt{-1}} 

\hyphenation{mul-ti-pli-ci-ty}

\maketitle

\begin{abstract}
Let $X$ be a compact Riemann surface together with a finite set of marked
points. We use Morse theoretic techniques to compute the Betti
numbers of the parabolic $U(2,1)$-Higgs bundles moduli spaces over
$X$. We give examples for one marked point showing that the
Poincar\'e polynomials depend on the system of weights of the
parabolic bundle.
\end{abstract}

\section{Introduction}

The moduli spaces of stable parabolic Higgs bundles have been studied
in \cite{by,t,y} and have a rich structure, partially due to its relation with
the representation space of the fundamental group of a punctured
Riemann surface. This relationship was established by Simpson in
\cite{s}.

The topology of the moduli $\cU$ of stable $U(p,q)$-parabolic Higgs
bundles with fixed generic weights and degrees has already been
considered in \cite{glm}.
This moduli space is a submanifold of the moduli space $\cM$ of stable parabolic
Higgs bundles of fixed degree, which was analised in the rank $2$ case
by Boden and Yokogawa in \cite{by} and in the rank $3$ case by
Garc{\iacute}a-Prada, Mu\~{n}oz and Gothen in \cite{ggm}. 
In these two papers, the authors obtained the Betti numbers of the
moduli spaces $\cM$. 
Our purpose here is to calculate Betti numbers of $\cU$ when $p+q=3$.
Note that in this case, $\cU$ is a submanifold of the moduli $\cM$
studied in \cite{ggm}. 
It is known that for fixed rank, the moduli spaces $\cM$ of stable
parabolic Higgs bundles, corresponding to different choices of
degrees and generic weights, are diffeomorphic \cite{ggm, t}.
Our computations for $p+q=3$ produce counterexamples to this type of
phenomena for the submanifolds $\cU$, 
that is, they provide an example of the dependence of these moduli
spaces on the generic weights of the parabolic structure.

We will use Morse theory one step forward than in \cite{glm} thanks
to fixing the rank equal to $3$. Higher ranks need to develop
another tool called parabolic chains and will be done in the future.
We start in Section \ref{sec:def}, explaining the necessary
definitions and results for the Morse theory involved and
defining also the Morse function that we are going to use. In Section
\ref{sec:111} we study certain critical subvarieties of this Morse
function before and in Section \ref{sec:2112} we introduce
parabolic triples for another type of critical subvarieties.
Sections \ref{sec:111-1} and \ref{12-1} give explicit computations
for the case of one puncture and Section \ref{sec:betti} summarizes
the results and give some low genus examples.

\subsection*{Acknowledgements} We wish to thank Vicente Mu\~{n}oz for
very useful comments and corrections, and to Luis \'Alvarez-Consul for 
his help.

\section{Definitions and Morse theory}\label{sec:def}

Let $X$ be a compact Riemann surface of genus $g\ge 0$ together with
a finite set of marked distinct
points $x_{1}, \ldots, x_{s}$. We denote $D=x_{1}+\cdots+x_{s}$ the
divisor on $X$ defined by the punctures.

A parabolic bundle $E$ over $X$ consists of a holomorphic bundle
with a parabolic structure, that is, weighted flags, one for each
puncture in $X$,
\begin{eqnarray*}
E_{x}&=&E_{x,1}\supset \cdots \supset E_{x,r(x)}\supset 0,\\
0&\le&\alpha_{1}(x)<\cdots < \alpha_{r(x)}<1.
\end{eqnarray*}
The set of all weights for all $x\in D$,
$\alpha=\{\alpha_{i}(x);i=1, \ldots, r(x)\}$, is called
\emph{parabolic system of weights} of $E$.

A holomorphic map $f:E\to E'$ between parabolic bundles is called
\emph{parabolic} if $\alpha_{i}(x)>\alpha'_{j}(x)$ implies
$f(E_{x,i})\subset E'_{x,j+1}$  for all $x\in D$, and  $f$
\emph{strongly parabolic} if $\alpha_{i}(x)\ge \alpha'_{j}(x)$
implies $f(E_{x,i})\subset E'_{x,j+1}$ for all $x\in D$, where we
denote by $\alpha'_{j}(x)$ the weights on $E'$. Also $\PH(E,E')$ and
$\SPH(E,E')$ will denote respectively the bundles of parabolic and
strongly parabolic morphisms from $E$ to $E'$. Finally, a
\emph{parabolic subbundle} of a parabolic bundle is a subbundle which
inherits its parabolic structure from the parabolic bundle.

We write $m_{\alpha_{i}}(x)=\dim (E_{x,i}/E_{x,i+1})$ for the
multiplicity of the weight $\alpha_{i}(x)$ at $x$.
The parabolic degree and parabolic slope of $E$ are defined as
$$\pdeg(E)= \deg(E)+\sum_{x\in
D}\sum_{i=1}^{r(x)}m_{\alpha_{i}}(x)\alpha_{i}(x),$$
$$\pmu(E)=\frac{\pdeg(E)}{\rk(E)}.$$

A parabolic bundle is called (semi)-stable if for every parabolic
subbundle $F$ of $E$, the parabolic slope satisfies $\pmu(F)\le
\pmu(E)$ (resp. $\pmu(F)<\pmu(E)$).

For parabolic bundles $E$ there is a well-defined notion of
parabolic dual $E^{\ast}$. It consists of the bundle $\Hom(E,\cO(-D))$ and at each $x\in D$ a weighted
filtration 
\begin{eqnarray*}
E^{\ast}_{x}=E^{\ast}_{x,1}\supset\cdots E^{\ast}_{x,r(x)}\supset
0,\\
0<1-\alpha_{r(x)}(x)<\cdots<1-\alpha_{1}<1.
\end{eqnarray*}
In the case $\alpha_{1}=0$ we choose
the following weights for the filtration,
$$0\le\alpha_{1}<1-\alpha_{r(x)}(x)<\cdots<1-\alpha_{2}<1.$$

With this definition $E^{\ast\ast}=E$ and
$\pdeg(E^{\ast})=-\pdeg(E)$.

A \emph{$GL(n,\CC)$-parabolic Higgs bundle} is a pair $(E,\Phi)$
consisting of a parabolic bundle $E$ and $\Phi\in
H^{0}(\SPE(E)\otimes K(D))$, i.e. $\Phi$ is a meromorphic
endomorphism valued one-form with simple poles along $D$ whose
residue at $p\in D$ is nilpotent with resect to the flag. A
parabolic Higgs bundle is called (semi)-stable if for every
$\Phi$-invariant subbundle $F$ of $E$, its parabolic slope
satisfies $\pmu(F)\le \pmu(E)$ (resp. $\pmu(F)<\pmu(E)$). We shall
say that the weights are \emph{generic} when every semistable
Higgs bundle is stable, that is, there are no properly semistable
parabolic Higgs bundles.

A \emph{$U(p,q)$-parabolic Higgs bundle} on $X$ is a parabolic
Higgs bundle $(E,\Phi)$ such that $E=V\oplus W$, where $V$ and $W$
are parabolic vector bundles of rank $p$ and $q$ respectively, and
  $$
  \Phi=\left(\begin{array}{ll}  0 & \beta \\ \gamma & 0
  \end{array}\right) :(V\oplus W)\to (V\oplus W)\otimes K(D),
  $$
where the non-zero components $\beta:W\to V\otimes K(D)$ and
$\gamma: V\to W\otimes K(D)$ are strongly parabolic morphisms.
Hence a $U(p,q)$-parabolic Higgs bundle is (semi)-stable if the
slope (semi)-stability condition is satisfied for all
$\Phi$-invariant subbundles of the form $F=V'\oplus W'$, i.e. for
all subbundles $V'\subset V$ and $W'\subset W$ such that
\begin{eqnarray}
\beta:&W'\to V'\otimes K(D)\\
\gamma:&V'\to W'\otimes K(D).
\end{eqnarray}

Let us fix generic weights and topological invariants $\rk(E)$ and
$\deg(E)$. The moduli space $\cM_{GL(n,\CC)}$ of stable
$GL(n,\CC)$-parabolic Higgs bundles was constructed using Geometric
Invariant Theory by Yokogawa \cite{y}, who also showed that it is a
smooth irreducible complex variety.

By definition there is an injection from the moduli $\cU_{(p,q)}$ of
stable $U(p,q)$-parabolic Higgs bundles to the moduli
$\cM_{GL(p+q,\CC)}$ of stable $GL(p+q,\CC)$-parabolic Higgs bundles.
Moreover, such an injection is an embedding, as shown in \cite{glm}, so
$\cU_{(p,q)}$ is in fact a submanifold of $\cM_{GL(p+q,\CC)}$. When
it does not induce confusion, we will denote $\cU_{(2,1)}$ and
$\cM_{GL(3,\CC)}$ by $\cU$ and $\cM$.

The Toledo invariant for the moduli of $U(p,q)$ parabolic Higgs
bundles is studied in \cite{glm} and defined as
$\tau=2(q\pdeg(V)-p\pdeg(W))/(p+q)$. Thus, given $(E,\Phi)\in \cU$
we have
\begin{equation}\label{def:toledo}
\tau=\frac{2}{3}(\Delta-3b+\sum_{x\in
D}\alpha_{1}(x)+\alpha_{2}(x)-2\eta(x)),
\end{equation}
where we denote $a=\deg(V)$, $b=\deg(W)$, $\alpha_{1}(x)$ and
$\alpha_{2}(x)$ the parabolic weights on $V$  and $\eta(x)$ the
parabolic weights on $W$ over the punctures $x\in D$, and 
$\Delta=a+b$. We will use this notation in the following.

\begin{prop}\label{prop:isos}
The map $V\oplus W\to (V\oplus W)\otimes L$, where $L$ is a
parabolic line bundle, induces
an isomorphism from the moduli space $\cU_{(p,q)}(a,b)$ of
parabolic $U(p,q)$-Higgs bundles with fixed degrees $(a,b)$ to the
moduli space $\cU_{(p,q)}(a',b')$ of parabolic $U(p,q)$-Higgs
bundles with fixed degrees $(a',b')$, where $a'=a+pl$ and $b'=b+ql$.

The map $V\oplus W\to V^{\ast}\oplus W^{\ast}$ induces 
an isomorphism of moduli spaces, from $\cU_{(p,q)}(a,b)$ to
$\cU'_{(p,q)}(a',b')$, where $a'=-a$ and $b'=-b$.  \hfill $\Box$
\end{prop}

Let us to assume that $\Delta=a'+b'\equiv 0(3)$ and that the parabolic Toledo
invariant $\tau$ satisfies $\tau\geq 0$.

The moduli space $\cU$ of stable $U(p,q)$-parabolic Higgs bundles
has been studied in \cite{glm}, where the number of connected
components is calculated using Bott-Morse theory. Here we shall fix
later $p=2$ and $q=1$ to go one step further and give topological
information about this moduli space.

Consider the action of $\CC^{\ast}$ on $\cU$ given in \cite{glm} as
\begin{eqnarray}
\psi: \CC^{\ast}\times \cU&\to& \cU\\
(\lambda,(E,\Phi))&\mapsto&  (E,\lambda\Phi).
\end{eqnarray}
This restricts to a Hamiltonian action of $S^{1}\subset \CC^{\ast}$
on $\cU$ and the moment map associated to this Hamiltonian action is
defined by
\begin{equation}\label{eq:f}
f([E,\Phi])=\|\Phi\|^{2}=\frac{1}{\pi}\|\beta\|^{2}+\frac{1}{\pi}\|\gamma\|^{2},
\end{equation}
where we are using a suitable Sobolev metric for the norm given by the Hermite-Einstein equations for the parabolic Higgs bundle $(E,\Phi)$ (see \cite{s}).

Observe that $f:\cU\to \RR$ is the restriction of the moment map $f:\cM\to \RR$ used in \cite{ggm}. That map was proper.  Hence, $f$ is also proper since $\cU$ is a closed submanifold of $\cM$. This fact together with a result of Frankel \cite{f}, proving that a proper moment map for a Hamiltonian circle action on a K\"{a}hler manifold is a perfect Bott-Morse function,  give us that $f$ is a perfect  Bott-Morse function. 

Hence, we have the following formula for the Poincar\'e polynomial
of the manifold $\cU$,
\begin{equation}\label{eq:Poincare}
P_{t}(\cU)=\sum_{\cN}t^{\lambda_{\cN}}P_{t}(\cN),
\end{equation}
where the sum runs over all critical submanifolds $\cN$ of $\cU$
for $f$ and $\lambda_{\cN}$ is the Morse index of $f$ on $\cN$.

The critical points of $f$ are exactly the fixed points of the
circle action. Moreover, the Morse index of $f$ at a critical point
equals the dimension of the negative weight space of the circle
action on the tangent space \cite{f}.

Simpson's theorem gives us a criterion for $(E,\Phi)$ to be a
critical point for the Morse function.

\begin{thm}[\cite{s}, Thm.8]\label{prop:simpson}
The equivalence class of a stable parabolic Higgs bundle
$(E,\Phi)$ is fixed under the action of $S^{1}$ if and only if it
is a \emph{parabolic complex variation of Hodge structure}. This
means that $E$ has a direct sum decomposition
$$E=E_{0}\oplus E_{1}\oplus \cdots \oplus E_{m}$$
as parabolic bundles, such that $\Phi$ is strongly parabolic and of
degree one with respect to this decomposition, in other words, the
restriction $\Phi_{l}=\Phi|_{E_{l}} \in
H^{0}(\SPH(E_{l},E_{l+1})\otimes K(D)).$ Also $\Phi_{l}\ne 0$ and
the weight of $\psi$ on $E_{l+1}$ is one plus the weight of $\psi$
on $E_{l}$.
\end{thm}

Finally, the Morse index of $f$ is calculated using the following
result.
\begin{prop}[\cite{glm}]\label{prop:complex-0}
The dimension of the eigenspace of the action of $\psi$ on the
tangent space for the eigenvalue $-k$ equals the first
hypercohomology group of a complex
$$C^{\bullet}_{k}:
U_{k}\to \bar{U}_{k}\otimes K(D)$$ where
\begin{equation}
\begin{array}{ll}
    U_{k}=\oplus_{j-i=2k}\PH(E_{i},E_{j}) &
    \bar{U}_{k}=\oplus_{j-i=2k+1}\SPH(E_{i},E_{j})
    \end{array}
\end{equation}
\end{prop}

Thus in the $U(2,1)$ case we have the following possibilities. If
$(E, \Phi)$ is a critical point then it can be of one of these
three following forms:
$$\begin{array}{l l}
                        E=E_{0}\oplus E_{1}, & \rk(E_{0})=1, \quad
                        \rk(E_{1})=2\\
                        E=E_{0}\oplus E_{1}, & \rk(E_{0})=2, \quad
                        \rk(E_{1})=1\\
                        E=E_{0}\oplus E_{1}\oplus E_{2} &
                        \rk(E_{i})=1, \quad  i=1,2,3.\end{array}$$
These form, critical subvarieties of types $(\rk(E_{0}),\rk(E_{1})$
or $(\rk(E_{0}),\rk(E_{1}),\rk(E_{2}))$ particularly in this case
$(1,2)$, $(2,1)$, and $(1,1,1)$ respectively. And this critical
subvarieties can be identified with triples of type
$(1,2,d_{1},d_{0};\alpha_{1},\alpha_{2},\eta)$,
$(2,1,d_{1},d_{0};\eta,\alpha_{1},\alpha_{2})$ and chains of type
$(1,1,1,d_{2},d_{1},d_{0};\alpha_{\varpi(1)},\eta,\alpha_{\varpi(2)})$.
Where by \emph{type} of a triple or a chain we mean, a system of
numbers that give some topological invariants of this objects, they
are the ranks, degrees and parabolic systems of weights of each
parabolic bundle conforming the triples or the chain respectively.

Observe that the critical varieties of type $(1,2)$ and type $(2,1)$
consist of parabolic Higgs bundles for which either $\gamma=0$ or
$\beta=0$ respectively. From (\ref{eq:f}) and using the definition
of $\tau$ we get that they are minima for the Morse function $f$
and, as proved in \cite{glm}, they are the only ones. Hence its Morse
index is zero.

In the cases $(2,1)$ and $(1,2)$ where $E=E_{0}\oplus E_{1}$ the
critical submanifold will be identified with certain moduli spaces
of parabolic triples. However in the third case, where
$E=E_{0}\oplus E_{1}\oplus E_{2}$, we will be dealing with parabolic
chains. This is the reason for restricting attention to $p=2$ and
$q=1$. If we would like to compute the Betti numbers for higher
values of $p$ and $q$ we will have to deal with more general
parabolic chains that the ones appearing here, and this tool has not
been developed yet. This is left to future work.

In the following sections we will calculate the Poincar\'e
polynomials which take part in the formula in (\ref{eq:Poincare}),
that is for the moduli space $\cU$ of parabolic $U(2,1)$-parabolic
Higgs bundles
\begin{equation}\label{eq:Poincare2}
P_{t}(\cU)=\left\{\begin{array}{lc}
            P_{t}\cN_{(2,1)}+
            P_{t}\cN_{(1,1,1)}&\mathrm{for\;}
            \tau>0\\
            P_{t}\cN_{(1,2)}+
            P_{t}\cN_{(1,1,1)}&\mathrm{for\;}
            \tau<0
            \end{array}\right.
\end{equation}
where we denote $P_{t}\cN_{(1,2)}$ the contribution on the
Poincar\'e polynomial of $\cU$ of the subvariety of type $(1,2)$,
$P_{t}\cN_{(2,1)}$ is the contribution of the subvariety of type
$(2,1)$ and $P_{t}\cN_{(1,1,1)}$ is the contributions from all
critical subvarieties of type $(1,1,1)$. Through these sections our
computations will depend on some variables that we have mentioned
above: the Toledo invariant $\tau$ of the moduli space, and
 the degree $\Delta=a+b$ of $E$. Recall that by Proposition
\ref{prop:isos} we can suppose $\tau<0$ and $\Delta\equiv 0(3)$.

It is known that for fixed rank, and for different choices of
degrees and generic weights the moduli spaces of parabolic Higgs
bundles $\cM$ have the same Poincar\'e polynomial (see \cite{ggm}),
so it is possible to choose the weights conveniently for such
calculation for $\cM$. We have seen that $\cU\subset \cM$ is a
subvariety, and our calculation of its Poincar\'e polynomial will
show that the same phenomenon does not happen for $\cU$. The
Poincar\'e polynomial of $\cU$ depends on the generic weights. We
shall see this very explicitly in our calculations for one marked
point.

\section{Contribution to Poincar\'e polinomial from critical subvarieties of type
(1,1,1).}\label{sec:111}
 We start with the case where $E=V\oplus W$ splits
in three line bundles $E=E_{0}\oplus E_{1} \oplus E_{2}$ where
$E_{0}$ and $E_{2}$ are contained in $V$, together with strongly
parabolic homomorphisms $\Phi_{0}=\gamma|_{E_{0}}:E_{0}\to
E_{1}\otimes K(D)$ and $\Phi_{1}=\beta|_{E_{1}}:E_{1}\to
E_{2}\otimes K(D)$.

We denote along this section $d_{i}=\deg(E_{i})$ so
$\Delta=d_{0}+d_{1}+d_{2}=d_{0}+b+d_{2}$ i.e. $a=d_{0}+d_{2}$ and
$b=d_{1}$.

The distributions of the weights for $E_{0}$ and $E_{2}$ are given
by a set of injective maps
$\varpi=\{\varpi_{x}:\{1,2\}\to\{1,2\};\;x\in D\}$ such that the
weight of $E_{0}$ at $x\in D$ is $\alpha_{\varpi(1)_{x}}(x)$ and the
weight of $E_{2}$ at $x\in D$ is $\alpha_{\varpi(2)_{x}}(x)$.

\begin{prop}The Morse index for the critical submanifolds of type (1,1,1) depends
on $d_{0}$ and $\varpi$, and it is given by
\begin{equation}
\lambda_{\cN_{(1,1,1)}}(d_{0},\varpi)=2g-2+2(2d_{0}-\Delta+b)+2(s-v))\\
\end{equation}
where $v=\sharp\{x\in
D;\;\alpha_{\varpi_{x}(1)}(x)>\alpha_{\varpi_{x}(2)}(x)\}$, and $s$
is the number of marked points.
\end{prop}

\begin{proof}
By proposition \ref{prop:complex-0} the Morse index equals the
dimension of $\HH^{1}(C^{\bullet}_{1})$ where $C^{\bullet}_{1}$ is
the complex
$$\PH(E_{0},E_{2})\to 0.$$

Using the long exact sequence for this complex we get
$\HH^{0}(C^{\bullet}_{1})=0$ since it is isomorphic to
$H^{0}(\PH(E_{0},E_{2}))$ and, the last is equal to zero since its
degree is less than zero. Hence,
\begin{eqnarray*}
\frac{1}{2}\lambda_{\cN_{(1,1,1)}}&=&\dim T_{E}\cU_{< 0}=\dim\HH^{1}(C^{\bullet}_{1})\\
&=&\dim
H^{1}(\PH(E_{0},E_{2}))=-\chi(\PH(E_{0},E_{2}))\\&=&-\deg(\PH(E_{0},E_{2}))-\rk(\PH(E_{0},E_{2}))(1-g)\\
&=& d_{0}-d_{2}+s-\sum_{x\in D}\dim\PH(E_{0},E_{2})_{x}+g-1.
\end{eqnarray*}

Hence, $\lambda_{\cN_{(1,1,1)}}=2g-2+2(2d_{0}+b-\Delta)+2(s-v)$,
where $v=\sharp\{x\in D; \alpha_{\varpi_{x}(1)}(x)\le
\alpha_{\varpi_{x}(2)}(x)\}$.
\end{proof}

\begin{rmk}The Proposition above  proves also that $\lambda_{\cN(1,1,1)}$ depend only
on $d_{0}$ and $\varpi$, the data that give us how splits $V$ into
$E_{0}$ and $E_{2}$. So we may decompose
$\cN(1,1,1)=\bigcup_{d_{0},\varpi}\cN(d_{0},\varpi)$.
\end{rmk}

From now on we denote
$$v_{1}=\sharp\{x\in
D;\,\alpha_{\varpi(1)}(x)<\eta(x)\}$$
$$v_{2}=\sharp\{x\in
D;\,\eta(x)<\alpha_{\varpi(2)}(x)\}.$$

\begin{prop}\label{prop:iso1}
Assume $\tau<0$, let $\cN_{(1,1,1)}$ be the union of critical
submanifolds of type $(1,1,1)$ parametrized by $d_{0}$ and $\varpi$,
i.e. $\cN_{(1,1,1)}=\bigcup_{d_{0},\varpi}\cN(d_{0},\varpi)$. The
map
\begin{eqnarray*}
\cN(d_0,\varpi)&\to& \Jac^{d_{0}}X\times S^{m_{1}}X\times
S^{m_{2}}X\\
(E_{0}\oplus E_{1}\oplus E_{2}, \Phi_{0},\Phi_{1})&\mapsto&(E_{0},
\di(\Phi_{0}),\di(\Phi_{1}))
\end{eqnarray*}
where
\begin{eqnarray*}
m_{1}&=&\deg(\SPH(E_{0},E_{1})\otimes
K(D))=b-d_{0}+2g-2+v_{1}\\
m_{2}&=&\deg(\SPH(E_{1},E_{2})\otimes
K(D))=\Delta-d_{0}-2b+2g-2+v_{2}
\end{eqnarray*}
is an isomorphism, in particular there is only one component for
fixed $d_{0}$ and $\varpi$. Furthermore, $d_{0}$ the degree of
$E_{0}$ is lower bounded by $\bar{d_{0}}$, that is,
\begin{equation}\label{eq:bsd0}
d_{0}\ge\bar{d_{0}}=\left[\frac{1}{3}\left(\Delta+\sum_{x\in D}
(\eta(x)+\alpha_{\varpi_{x}(2)}(x)-2\alpha_{\varpi_{x}(1)}(x))\right)+1\right]
\end{equation}
where $[k]$ denote the entire part of $k$.
\end{prop}

\begin{proof}

The isomorphism is obvious (see \cite{ggm}). The stability condition
on $E$ applied on the subbundles $E_{2}$ and $E_{1}\oplus E_{2}$,
together with the formula $d_{2}=\Delta-b-d_{0}$ gives the following
two bounds for $d_{0}$:
\begin{eqnarray}\label{eq:bd0}
2\Delta-3b-\sum_{x\in
D}(\alpha_{\varpi_{x}(1)}(x)+\eta(x)-2\alpha_{\varpi_{x}(2)}(x))<3d_{0}\\
\Delta-\sum_{x\in D}
(2\alpha_{\varpi_{x}(1)}(x)-\eta(x)-\alpha_{\varpi_{x}(2)}(x))<3d_{0}.
\end{eqnarray}
To determine which is the appropriate bound we subtract these two
inequalities. This subtraction gives a multiple of $\tau$ , hence
$\bar{d}_{0}$ depends on wether $\tau$ is negative or positive.
\end{proof}

\begin{rmk} The condition on the weights being generic implies that
$\tau$ can not be zero. This is because $\tau=0$ implies that
$2\eta(x)-\alpha_{1}(x)-\alpha_{2}(x)=\Delta-3b$, and if that
happens then there is a $U(2,1)$-parabolic Higgs  subbundle $(V,
\Phi=0)$ non-stable but semistable.
\end{rmk}
\begin{rmk}
Note that the values $m_{1}$ and $m_{2}$ depend on $\varpi_{x}$ and
$\bar{d_{0}}$.
\end{rmk}

\begin{rmk}We chose $\tau<0$ for computability reasons.
\end{rmk}

Denote $\varpi=\{\varpi_{x}\}_{x\in D}$.

\begin{thm}\label{cor:P-t-N-111}

  The Poincar\'e polynomial of the critical submanifold
  $\cN(d_0,\varpi)$ is
  \begin{displaymath}
    P_{t}(\cN(d_0,\varpi)) = (1+t)^{2g}
    \Coeff_{x^{0}y^{0}}\left(
    \frac{(1+xt)^{2g}}{(1-x)(1-xt^2)x^{m_{1}}}
    \cdot
    \frac{(1+yt)^{2g}}{(1-y)(1-yt^2)y^{m_{2}}}
    \right)
  \end{displaymath}
  where $m_{1}$ and $m_{2}$ are the same as in Proposition
  \ref{prop:iso1}.
\end{thm}

\begin{proof}
Use  Macdonald's formula for the Poincar\'e polynomial of the
symetric product (see \cite{md}).
\end{proof}

Now, in order to get the contribution of all the subvarieties of
type $(1,1,1)$ in $P_{t}(\cU)$ we have to sum over all $d_{0}\ge
\bar{d_{0}}$ and all possibilities of $\varpi$.

\begin{equation*}\begin{split}
&P_{t}(\cN_{(1,1,1)})=\sum_{d_{0},\varpi}t^{\lambda_{\cN(d_{0},\varpi)}}P_{t}(\cN(d_{0},\varpi))\\
&=\sum_{d_{0},\varpi}\left(t^{2g-2+2(b-\Delta)+4d_{0}+2(s-v)}
\Coeff_{x^{0}y^{0}}\left(
    \frac{(1+xt)^{2g}}{(1-x)(1-xt^2)x^{m_{1}}}
    \cdot
    \frac{(1+yt)^{2g}}{(1-y)(1-yt^2)y^{m_{2}}}\right)\right)\\
&=\Coeff_{x^{0}y^{0}}\left(\sum_{\varpi}\frac{t^{2g-2+2b-2\Delta+2s}(1+xt)^{2g}(1+yt)^{2g}}
{(1-x)(1-xt^{2})x^{b+2g-2}(1-y)(1-yt^2)y^{\Delta-2b+2g-2}} \cdot
\frac{t^{4\bar{d_{0}}}x^{\bar{d_{0}}}y^{\bar{d_{0}}}}{t^{2v}x^{v1}y^{v2}}\right)\\
&=\Coeff_{x^{0}y^{0}}\left(\frac{t^{2g-2+2b-2\Delta+2s}(1+xt)^{2g}(1+yt)^{2g}}
{(1-x)(1-xt^{2})x^{b+2g-2}(1-y)(1-yt^2)y^{\Delta-2b+2g-2}} \cdot
\sum_{\varpi}\frac{t^{4\bar{d_{0}}}x^{\bar{d_{0}}}y^{\bar{d_{0}}}}{t^{2v}x^{v1}y^{v2}}\right)
\end{split}
\end{equation*}

Thus, we have to compute the following sum
\begin{equation}\label{eq:suma}
\sum_{\varpi_{x}}\frac{t^{4\bar{d_{0}}}x^{\bar{d_{0}}}y^{\bar{d_{0}}}}{t^{2v}x^{v_{1}}y^{v_{2}}}.
\end{equation}
The variables depend also on the weights $\alpha_{1}(x)$,
$\alpha_{2}(x)$ and $\eta(x)$, and the distribution functions
$\varpi_{x}$.

\section{Computations for one puncture for
$\cN_{(1,1,1)}$.}\label{sec:111-1} From now on we consider the case
of one puncture to get more explicit formulas, so we denote
$\alpha_{i}=\alpha_{i}(x)$ for $i=1,2$ and $\eta=\eta(x)$. We
abbreviate $\varpi_{x}$ to $\varpi$.

We have to consider the following cases for the possible
distributions of the weights,

\begin{table}[h]
  \centering
  \caption{Weight distributions.}
\begin{tabular}{|c|c|}
\hline $S_{1}$& $\eta<\alpha_{1}<\alpha_{2}$\\
$S_{1}(a)$ & $\alpha_{2}-\alpha_{1}>\alpha_{1}-\eta$\\
$S_{1}(b)$ & $\alpha_{2}-\alpha_{1}<\alpha_{1}-\eta$\\
\hline $S_{3}$ & $\alpha_{1}<\alpha_{2}<\eta$\\
$S_{3}(a)$ & $\alpha_{2}-\eta>\alpha_{2}-\alpha_{1}$\\
$S_{3}(b)$ & $\alpha_{2}-\eta<\alpha_{2}-\alpha_{1}$\\
\hline
\end{tabular}
\qquad
\begin{tabular}{|c|c|}
\hline $S_{2}$ & $\alpha_{1}<\eta<\alpha_{2}$\\
\hline $S_{4}$& $\eta=\alpha_{1}<\alpha_{2}$\\
\hline $S_{5}$ & $\alpha_{1}<\alpha_{2}=\eta$\\
\hline $S_{6}$ & $\alpha_{1}=\alpha_{2}=\eta$\\
\hline $S_{7}$ &$\eta<\alpha_{1}=\alpha_{2}$\\
\hline $S_{8}$ & $\alpha_{1}=\alpha_{2}< \eta$\\
\hline
\end{tabular}
\label{tab:distr}
\end{table}

\begin{thm} The contributions to the Poincar\'e polynomial of the union of the subvarieties of type
$(1,1,1)$ when $\tau<0$ and one marked point are classified by the
possibilities for the distribution of the weights of $E$ shown in
Table \ref{tab:distr}. They are the following,

\begin{itemize}
\item[(i)] For $S_{1}(a)$:
$$\Coeff_{x^{0}y^{0}} \left(\frac{t^{2\,b - \frac{2\,\Delta}{3} + 2\,g}\,
    x^{2 - b + \frac{\Delta}{3} - 2\,g}\,
    {\left( 1 + t\,x \right) }^{2\,g}\,
    y^{1 + 2\,b - \frac{2\,\Delta}{3} - 2\,g}\,
    {\left( 1 + t\,y \right) }^{2\,g}\,
    \left( 1 + t^2\,x\,y \right) }{\left( -1 + x \right) \,
    \left( -1 + t^2\,x \right) \,\left( -1 + y \right) \,
    \left( -1 + t^2\,y \right) }\right). $$

\item[(ii)] For $S_{1}(b)$ :
$$\Coeff_{x^{0}y^{0}} \left(\frac{t^{-2 + 2\,b - \frac{2\,\Delta}{3} + 2\,g}\,
    \left( 1 + t^2 \right) \,
    x^{2 - b + \frac{\Delta}{3} - 2\,g}\,
    {\left( 1 + t\,x \right) }^{2\,g}\,
    y^{1 + 2\,b - \frac{2\,\Delta}{3} - 2\,g}\,
    {\left( 1 + t\,y \right) }^{2\,g}}{\left( -1 + x \right)
      \,\left( -1 + t^2\,x \right) \,
    \left( -1 + y \right) \,\left( -1 + t^2\,y \right) }\right). $$

\item[(iii)] For $S_{2}$ :
$$\Coeff_{x^{0}y^{0}} \left(\frac{t^{2\,b - \frac{2\,\Delta}{3} + 2\,g}\,
    \left( 1 + t^2 \right) \,
    x^{2 - b + \frac{\Delta}{3} - 2\,g}\,
    {\left( 1 + t\,x \right) }^{2\,g}\,
    y^{2 + 2\,b - \frac{2\,\Delta}{3} - 2\,g}\,
    {\left( 1 + t\,y \right) }^{2\,g}}{\left( -1 + x \right)
      \,\left( -1 + t^2\,x \right) \,
    \left( -1 + y \right) \,\left( -1 + t^2\,y \right) }\right).$$

\item[(iv)] For $S_{3}(a)$ :
$$\Coeff_{x^{0}y^{0}} \left(\frac{t^{2 + 2\,b - \frac{2\,\Delta}{3} + 2\,g}\,
    \left( 1 + t^2 \right) \,
    x^{2 - b + \frac{\Delta}{3} - 2\,g}\,
    {\left( 1 + t\,x \right) }^{2\,g}\,
    y^{3 + 2\,b - \frac{2\,\Delta}{3} - 2\,g}\,
    {\left( 1 + t\,y \right) }^{2\,g}}{\left( -1 + x \right)
      \,\left( -1 + t^2\,x \right) \,
    \left( -1 + y \right) \,\left( -1 + t^2\,y \right) }\right).$$

\item[(v)] For $S_{3}(b)$:
$$\Coeff_{x^{0}y^{0}} \left(\frac{t^{2\,b - \frac{2\,\Delta}{3} + 2\,g}\,
    x^{1 - b + \frac{\Delta}{3} - 2\,g}\,
    {\left( 1 + t\,x \right) }^{2\,g}\,
    y^{2 + 2\,b - \frac{2\,\Delta}{3} - 2\,g}\,
    {\left( 1 + t\,y \right) }^{2\,g}\,
    \left( 1 + t^2\,x\,y \right) }{\left( -1 + x \right) \,
    \left( -1 + t^2\,x \right) \,\left( -1 + y \right) \,
    \left( -1 + t^2\,y \right) }\right).$$

\item[(vi)] For $S_{4}$ :
$$\Coeff_{x^{0}y^{0}} \left(\frac{t^{2\,b - \frac{2\,\Delta}{3} + 2\,g}\,
    x^{2 - b + \frac{\Delta}{3} - 2\,g}\,
    {\left( 1 + t\,x \right) }^{2\,g}\,
    \left( 1 + t^2\,x \right) \,
    y^{2 + 2\,b - \frac{2\,\Delta}{3} - 2\,g}\,
    {\left( 1 + t\,y \right) }^{2\,g}}{\left( -1 + x \right)
      \,\left( -1 + t^2\,x \right) \,
    \left( -1 + y \right) \,\left( -1 + t^2\,y \right) }\right).$$

\item[(vii)] For $S_{5}$ :
$$\Coeff_{x^{0}y^{0}} \left(\frac{t^{2\,b - \frac{2\,\Delta}{3} + 2\,g}\,
    x^{2 - b + \frac{\Delta}{3} - 2\,g}\,
    {\left( 1 + t\,x \right) }^{2\,g}\,
    y^{2 + 2\,b - \frac{2\,\Delta}{3} - 2\,g}\,
    {\left( 1 + t\,y \right) }^{2\,g}\,
    \left( 1 + t^2\,y \right) }{\left( -1 + x \right) \,
    \left( -1 + t^2\,x \right) \,\left( -1 + y \right) \,
    \left( -1 + t^2\,y \right) }\right).$$

\item[(viii)] For $S_{6}$ :
$$\Coeff_{x^{0}y^{0}} \left(\frac{2\,t^{2 + 2\,b - \frac{2\,\Delta}{3} + 2\,g}\,
    x^{3 - b + \frac{\Delta}{3} - 2\,g}\,
    {\left( 1 + t\,x \right) }^{2\,g}\,
    y^{3 + 2\,b - \frac{2\,\Delta}{3} - 2\,g}\,
    {\left( 1 + t\,y \right) }^{2\,g}}{\left( -1 + x \right)
      \,\left( -1 + t^2\,x \right) \,
    \left( -1 + y \right) \,\left( -1 + t^2\,y \right) }\right).$$

\item[(ix)] For $S_{7}$ :
$$\Coeff_{x^0 y^0}\left(\frac{2\,t^{-2 + 2\,b - \frac{2\,\Delta}{3} + 2\,g}\,
    x^{2 - b + \frac{\Delta}{3} - 2\,g}\,
    {\left( 1 + t\,x \right) }^{2\,g}\,
    y^{1 + 2\,b - \frac{2\,\Delta}{3} - 2\,g}\,
    {\left( 1 + t\,y \right) }^{2\,g}}{\left( -1 + x \right)
      \,\left( -1 + t^2\,x \right) \,
    \left( -1 + y \right) \,\left( -1 + t^2\,y \right) }\right).$$

\item[(x)] For $S_{8}$ :
$$\Coeff_{x^0 y^0}\left(\frac{2\,t^{-2 + 2\,b +
       4\,\left( 1 + \frac{\Delta}{3} \right)  -
       2\,\Delta + 2\,g}\,
    x^{2 - b + \frac{\Delta}{3} - 2\,g}\,
    {\left( 1 + t\,x \right) }^{2\,g}\,
    y^{3 + 2\,b - \frac{2\,\Delta}{3} - 2\,g}\,
    {\left( 1 + t\,y \right) }^{2\,g}}{\left( 1 - x \right)
      \,\left( 1 - t^2\,x \right) \,\left( 1 - y \right) \,
    \left( 1 - t^2\,y \right) }  \right).$$

\end{itemize}
\end{thm}

\begin{proof}
Compute the values of $\bar{d_{0}}$, $v_{1}$. $v_{2}$ and $v$ for
each possible distribution of the weights, then we obtain the values
for the sum in (\ref{eq:suma}) for each case $S_{i}$.

The value for $\bar{d_{0}}$ depends on the distribution of the
weights and is, in the case $\varpi=\id$,
$\bar{d_{0}}=\frac{\Delta}{3}+1$ for all $S_{i}$ except for
$S_{1}(b)$ and $S_{7}$ where $\bar{d_{0}}=\frac{\Delta}{3}$. When
$\varpi\ne\id$, $\bar{d_{0}}=\frac{\Delta}{3}$ for all $S_{i}$
except for $S_{3}(a)$, $S_{6}$ and $S_{8}$ where it is
$\bar{d_{0}}=\frac{\Delta}{3}+1$.

\end{proof}

\section{Poincar\'e Polynomial for critical subvarieties of type $(1,2)$.}\label{sec:2112}

Following our previous discussion the critical subvarieties of type
$(1,2)$ and $(2,1)$ can be identified with the moduli of
$(2g-2)$-stable parabolic triples of type
$(2,1,a+4g-4,b;\alpha_{1},\alpha_{2},\eta)$ and
$(1,2,b+2g-2,a;\alpha_{1},\alpha_{2},\eta)$ respectively. So we
recall the basics of parabolic triples from \cite{ggm}.

>From Proposition \ref{prop:isos} we restrict to the case when
$\tau>0$, note that by definition the Morse function $f$ forces
$\gamma=0$ when $\tau>0$. Hence, for our analysis we only have to
consider the critical subvarieties of type $(1,2)$, that is
$(2g-2)$-stable parabolic triples of type
$(2,1,a+4g-4,b;\alpha_{1},\alpha_{2},\eta)$.

A \emph{parabolic triple} is a holomorphic triple $T=(T_{1},
T_{2},\phi)$ where $T_{1}$ and $T_{2}$ are parabolic bundles over
$X$, and $\phi:T_{2}\to T_{1}(D)$ is a strongly parabolic
homomorphism, i.e. an element $\phi\in H^{0}(\SPH(T_{2},T_{1}(D)))$.
We call \emph{parabolic system of weights} for the triple $(T,\phi)$
to the vector $\alpha=(\alpha^{1},\alpha^{2})$ where $\alpha^{i}$ is
the system of weights of $T_{i}$ with $i=1,2$. The \emph{type} of a
parabolic triple is a $n$-tuple
$(r_{1},r_{2},d_{1},d_{2};\alpha_{1}(x),\ldots, \alpha_{r(x)}(x),
\eta_{1}(x),\ldots, \eta_{r'(x)}(x)$, where $r_{i}=\rk(T_{i})$,
$d_{i}=\deg(T_{i})$, $\alpha$ is the parabolic system of weights of
$T_{1}$ and $\eta$ is the parabolic system of weights of $T_{2}$.

A parabolic triple $T'=(T'_{1}, T'_{2},\phi')$ is a
\emph{parabolic subtriple} of $T=(T_{1}, T_{2},\phi)$ if
$T'_{i}\subset T_{i}$ are parabolic subbundles for $i=1,2$ and
$\phi'(T'_{2})\subset T'_{1}(D)$ where $\phi'$ is the restriction
of $\phi$ to $T'_{2}$.

For any $\sigma\in \RR$ the \emph{$\sigma$-parabolic degree} of
$T$ is defined to be
$$\pdeg_{\sigma}(T)=\pdeg(T_{1})+\pdeg(T_{2})+\sigma \rk(T_{2}).$$

In the following we denote $r_{1}=\rk(T_{1})$ and
$r_{2}=\rk(T_{2})$. Thus we have a notion of stability for a fixed
parameter. Let $\sigma$ be a real number. We define the
\emph{$\sigma$-slope of a triple} $(T_{1},T_{2},\phi)$ as
\begin{equation}\label{eq:sigmamu}
\pmu_{\sigma}(T)=\frac{\pdeg T_{1}+\pdeg
T_{2}}{r_{1}+r_{2}}+\sigma\frac{r_{2}}{r_{1}+r_{2}}.
\end{equation}
$T$ is called $\sigma$-stable (resp. $\sigma$-semistable) if for
any non-zero proper subtriple $T'$ we have
 $\pmu_{\sigma}(T')<\pmu_{\sigma}(T)$ (resp.$\le$).

\begin{prop}
Subvarieties of type $(1,2)$ and type $(2,1)$ correspond with
$\sigma$-stable triples for $\sigma=2g-2$.
\end{prop}
\begin{proof}
In the case of the study of critical varieties of type $(1,2)$,
Simpson's theorem says that we have a variation of the Hodge
structure like this:
$$E=E_{0}\oplus E_{1},\quad \Phi=\gamma:E_{0}=W\to E_{1}\otimes K(D)=V\otimes K(D),$$
with $\rk(E_{0})=1$ and $\rk(E_{1})=2$. Therefore we get
$T=(E_{1}\otimes K, E_{0},\beta)$ a parabolic triple of type
$(2,1,a+4g-4,b;\alpha_{1},\alpha_{2},\eta)$.

Analogously in the case of the study of critical varieties of type
$(2,1)$ Simpson's theorem give us a variation of the Hodge
structure like before.
$$E_{0}\oplus E_{1}\quad \Phi=\beta:E_{0}=V\to E_{1}\otimes K(D)=W\otimes K(D)$$
With $\rk(E_{0})=2$ and $\rk(E_{1})=1$.

Hence in the case of critical varieties of type $(2,1)$ we have to
study parabolic triples of type $(1,2,
b+2g-2,a;\eta,\alpha_{1},\alpha_{2})$ with
$T=(T_{1},T_{2},\phi)=(E_{1}\otimes K, E_{0},\beta)$.
\end{proof}

\begin{prop}
The Morse index for critical submanifolds of type $(2,1)$ and type
$(1,2)$ is $\lambda_{\cN}=0$. In particular it does not depend on
the weights.
\end{prop}

\begin{proof}
This is clear since these subvarieties are minima for the Morse
function.
\end{proof}

\section{Computations for one puncture for
$\cN_{(1,2)}$.}\label{12-1}

Let $\varpi$ be a fixed distribution of the weights over the marked
point $x$. In the following $v_{1}$, $v_{2}$ and $v_{3}$ are given
by
\begin{eqnarray*}
v_{1}&=&\left\{\begin{array}{ll}
                        1 & \mathrm{if}\quad
                        \eta<\alpha_{\varpi(2)}\\
                        0 & \mathrm{otherwise}
                     \end{array}\right.\\
v_{2}&=&\left\{\begin{array}{ll}
                        1 & \mathrm{if}\quad
                        \eta<\alpha_{\varpi(1)}\\
                        0 & \mathrm{otherwise}
                     \end{array}\right.\\
v_{3}&=&\left\{\begin{array}{ll}
                        1 & \mathrm{if}\quad
                        \alpha_{\varpi(1)}<\alpha_{\varpi(2)}\\
                        0 & \mathrm{otherwise}
                     \end{array}\right.\\
\end{eqnarray*}

Let $\sigma>\sigma_{m}$ be a non-critical value. For any $\varpi$,
$\bar{d_{M}}$ we define,
\begin{equation*}
\bar{d}_{M}=\left[\frac{1}{3}\left(\Delta+\alpha_{\varpi(2)}+\eta-2\alpha_{\varpi(1)}+\sigma\right)+1\right].
\end{equation*}
\begin{prop}\label{eq:PPtriples}
The Poincar\'e polynomial of the moduli of parabolic triples
$T=(T_{1},T_{2},\phi)$ of type
$(2,1,\bar{d_{1}},\bar{d_{2}};\alpha,\eta)$ and one marked point is
\begin{equation}
\Coeff_{x^{0}}\frac{(1+t)^{4g}(1+xt)^{2g}}{(1-t^{2})(1-x)(1-xt^{2})}
\sum_{\varpi}x^{\bar{d}_{M}-\bar{d}_{1}+\bar{d}_{2}-v_{1}}
\left(\frac{t^{2\bar{d}_{1}-2\bar{d}_{2}+2v_{2}+2v_{3}-2\bar{d}_{M}}}{1-t^{-2}x}
-\frac{t^{-2\bar{d}_{1}+2g-2v_{3}+4\bar{d}_{M}}}{1-t^{4}x}\right)
\end{equation}
\end{prop}
\begin{proof} Rewrite theorem 6.5 from \cite{ggm} for this concrete
conditions.
\end{proof}

\begin{rmk}In \cite{ggm} the Poincar\'e polynomial is computed under the
assumption of generic distinct weights but this formula does not use
the assumption.
\end{rmk}

\begin{thm} The Poincar\'e polynomials for the
critical variety of type $(1,2)$ when $\tau>0$ for one marked point
are classified by the possibilities for the distribution of the
weights of $E$ given in Table \ref{tab:distr} and, they are the
following,
\small
\begin{itemize}
\item[(i)] For $S_{1}(a)$ and $S_{7}$
\begin{eqnarray*}
\Coeff_{x^{0}}\left(\frac{(1+t)^{4\,g}\,x^{1+2\,b-\frac{2}{3}\Delta-2\,g}(1+t\,x)^{2\,g}}
{t^{4\,b +\frac{2}{3}\,\Delta-2\,g}(-1+t^2)\,(t^2-x)\,(-1+x)\,(-1+t^2\,x)\,(-1+t^4\,x)}\right.\\
\left.\left(t^{6\,b}\,x-t^{4+6\,b}\,x+t^{2\,\Delta+2\,g}\,(1+x)-t^{4\,+2\,\Delta+2\,g}\,x\,(1+x)
+t^{2+6\,b}\,(-1+x^2)\right)\right).
\end{eqnarray*}

\item[(i)] For $S_{1}(b)$
$$\Coeff_{x^{0}}\left(- \frac{{\left( 1 + t \right) }^{4\,g}\,
      \left( 1 + t^2 \right) \,
      x^{1 + 2\,b - \frac{2}{3}\Delta - 2\,g}\,
      {\left( 1 + t\,x \right) }^{2\,g}\,
      \left( t^{2 + 6\,b} - t^{ 2 + 2\,\Delta +2\, g  } -
        t^{6\,b}\,x + t^{6\, + 2\,\Delta + 2\,g  }\,x
        \right) }{t^
       { 2 + 4\,b + \frac{2}{3}\Delta - 2\,g }\,
      \left( -1 + t^2 \right) \,\left( t^2 - x \right) \,
      \left( -1 + x \right) \,\left( -1 + t^2\,x \right) \,
      \left( -1 + t^4\,x \right) }
\right).$$

\item[(ii)] For $S_{2}$ and $S_{4}$
$$\Coeff_{x^{0}}\left(-\frac{{\left( 1 + t \right) }^{4\,g}\,
      \left( 1 + t^2 \right) \,
      x^{2 + 2\,b - \frac{2}{3}\Delta - 2\,g}\,
      {\left( 1 + t\,x \right) }^{2\,g}\,
      \left( t^{4 + 6\,b} - t^{2\, \Delta + 2\,g  } -
        t^{2 + 6\,b}\,x + t^{ 4 + 2\,\Delta + 2\,g  }\,x
        \right) }{t^
       { 2 + 4\,b + \frac{2}{3}\Delta - 2\,g }\,
      \left( -1 + t^2 \right) \,\left( t^2 - x \right) \,
      \left( -1 + x \right) \,\left( -1 + t^2\,x \right) \,
      \left( -1 + t^4\,x \right) } \right).$$

\item[(iii)] For $S_{3}(a)$, $S_{6}$ and $S_{8}$
$$\Coeff_{x^{0}}\left(- \frac{{\left( 1 + t \right) }^{4\,g}\,
      \left( 1 + t^2 \right) \,
      x^{3 + 2\,b - \frac{2}{3}\Delta - 2\,g}\,
      {\left( 1 + t\,x \right) }^{2\,g}\,
      \left( t^{8 + 6\,b} - t^{2\,\Delta +2\, g  } -
        t^{6 + 6\,b}\,x + t^{4 + 2\,\Delta + 2\,g }\,x
        \right) }{t^
       {4 + 4\,b + \frac{2}{3}\Delta - 2\,g }\,
      \left( -1 + t^2 \right) \,\left( t^2 - x \right) \,
      \left( -1 + x \right) \,\left( -1 + t^2\,x \right) \,
      \left( -1 + t^4\,x \right) } \right).$$

\item[(iv)] For $S_{3}(b)$ and $S_{5}$
\begin{eqnarray*}
\Coeff_{x^{0}}\left(\frac{(1+\,t)^{4\,g}\,x^{2+2\,b-\frac{2\,\Delta}{3}-2\,g}\,(1+t\,x)^{2\,g}\,}
{t^{2+4\,b+\frac{2}{3}\Delta-2\,g}\,(-1+t^2)\,(t^2-x)\,(-1+x)\,(-1+t^2\,x)\,(-1+t^4\,x)}\right.\\
\left.\left(t^{2+6\,b}\,x-t^{6+6\,b}\,x+t^{2\,\Delta+2\,g}\,(1+x)-t^{4+2\,\Delta+2\,g}\,x\,(1+x)+t^{4+6\,b}\,(-1+x^2)\right)\right).
\end{eqnarray*}

\end{itemize}
\normalsize
\end{thm}

\begin{proof} We only have to apply Proposition \ref{eq:PPtriples} using the  different values for
$v_{1}$ and $v_{2}$  on each case. Use also that $v_{3}=1$ if
$\varpi=\id$ and equal to zero otherwise. Hence we have different
values for $\bar{d}_{M}$ depending on the distribution of the
weights. These are, when $\varpi=\id$,
$\bar{d_{M}}=\frac{\Delta}{3}+2g-1$ fora all $S_{i}$ except for
$S_{1}(a)$ and $S_{7}$ where $\bar{d_{M}}=\frac{\Delta}{3}+2g-2$.
And when $\varpi\ne\id$, $\bar{d_{M}}=\frac{\Delta}{3}+2g-2$ for all
$S_{i}$ except for $S_{6}$ and $S_{8}$ where
$\bar{d_{M}}=\frac{\Delta}{3}+2g-1$.
\end{proof}

\section{Poincar\'e polinomial of $\cU$ with one marked
point}\label{sec:betti}

Summarizing, we are using Morse-Bott theory in order to calculate
the Poincar\'e polynomial of $\cU$ the moduli of stable $U(2,1)$
parabolic Higgs bundles with degrees $a=\deg(V)$ and $b=\deg(W)$.
Therefore we have described the critical subvarieties of $\cU$ for
the Morse function $f$. These consist of  several subvarieties of
type $(1,1,1)$ parametrized by $(d_{0},\varpi)$ and, depending on
$\tau$, and one subvariety corresponding to the minima of $f$, which
is of type $(1,2)$ when $\tau<0$ and  of type $(2,1)$ when $\tau>0$.

\begin{cor}
The Poincar\'e polinomial of $U(2,1)$ parabolic  Higgs bundles for
$\tau<0$ is given by
$$P_{t}(\cU)=P_{t}(\cN_{(1,1,1)})+P_{t}(\cN_{2g-2}(2,1,a+4g-4,b))$$
\end{cor}

\begin{cor}
The Poincar\'e polynomial of $U(2,1)$ parabolic Higgs bundles for
$\tau>0$ is given by
$$P_{t}(\cM(a,b))=P_{t}(\cN_{(1,1,1)})+P_{t}(\cN_{2g-2}(1,2,b+2g-2,a))$$
\end{cor}

We can compute the Poincar\'e polynomial of the moduli space of
parabolic $U(2,1)$-Higgs bundles for specific values of $a$, $b$,
and $g$ using a computer algebra system.

Note that in order to give an example  we fix $a$, $b$ (so we fix
$\Delta$) such that $\Delta$ is equal to zero modulo $3$ and
$\tau<0$. The other varieties, with different values for $\Delta$
and $tau$ are diffeomorphic to these from Proposition
\ref{prop:isos}

Fix $g=1$, degrees $a=b=0$, and $\alpha_{1}<\alpha_{2}<\eta$ such
that $\eta-\alpha_{2}<\alpha_{2}-\alpha_{1}$, that is, we are in
case $S_{3}(b)$ of Table \ref{tab:distr}.

The contribution from the critical subvarieties of type $(1,1,1)$ is
$$
P_{t}(\cN_{(1,1,1)})=t^2 + 2\,t^3 + t^4,
$$
and from the critical subvariety of type $(1,2)$ is
$$
P_{t}(\cN_{2g-2}(2,1,a+4g-4,b))=1 + 4\,t + 6\,t^2 + 4\,t^3 + t^4.
$$

Hence, the Poincar\'e polynomial for $\cU(0,0)$ with $g=1$, when
$\alpha_{1}<\alpha_{2}<\eta$ such that
$\eta-\alpha_{2}<\alpha_{2}-\alpha_{1}$,
 is
$$
P_{t}(\cU)=1 + 4\,t + 7\,t^2 + 6\,t^3 + 2\,t^4.
$$

As an example of the phenomena we have talked above, different
polynomials for different weights we give another example for same
genus and degrees but for case $S_5$.

The contribution from the critical subvarieties of type $(1,1,1)$ is
$$
P_{t}(\cN_{(1,1,1)})=t^2,
$$
and from the critical subvariety of type $(1,2)$ is again
$$
P_{t}(\cN_{2g-2}(2,1,a+4g-4,b))=1 + 4\,t + 6\,t^2 + 4\,t^3 + t^4.
$$

Hence, the Poincar\'e polynomial for $\cU(0,0)$, when $g=1$, one
marked point and, $\alpha_{1}<\alpha_{2}= \eta$,  is
$$
P_{t}(\cU)=1+4\, t+7\, t^2+4\, t^3+t^4
$$

Note that can only choose degrees $a=b=0$ for
distributions of weights $S_{3}$ , $S_5$ and $S_8$

\begin{prop}
The complex dimension of the moduli space of parabolic
$U(2,1)$-Higgs bundles is $1+9(g-1)+\sum_{x\in D}(3-c)$ where $c$ is
the number of weights $\alpha_{i}(x)$ equal to $\eta(x)$.
\end{prop}
\begin{proof}
Rewrite Proposition $3.4$ from \cite{glm}.
\end{proof}

Hence, in the examples above the real dimension of $\cU$ is $6$ and
does not coincides with the degree of the polynomials. Also it is
interesting the fact that one of the polynomials satisfies
Poincar\'e duality and the other does not.

\end{document}